\documentclass[12pt]{amsart}
\usepackage{amsfonts, amssymb, amsmath, amsthm, latexsym, array}
\usepackage{fullpage,color,longtable}
\usepackage{verbatim,euscript}

\numberwithin{equation}{section}

\input {cyracc.def}
 \font\tencyr=wncyr10 
\font\tencyi=wncyi10 
\font\tencysc=wncysc10 
\def\rus{\tencyr\cyracc}
\def\rusi{\tencyi\cyracc}
\def\rusc{\tencysc\cyracc}

\newtheorem{thm}{Theorem}[section]
\newtheorem{prob}[thm]{Problem}
\newtheorem{lm}[thm]{Lemma}
\newtheorem{cl}[thm]{Corollary}
\newtheorem{prop}[thm]{Proposition}

\theoremstyle{remark}
\newtheorem{ex}[thm]{Example}
\newtheorem*{rem}{Remark}

\theoremstyle{definition}
\newtheorem{rmk}[thm]{Remark}
\newtheorem{df}{Definition}

\newenvironment{E6}[6]{%
{\small\begin{tabular}{@{}c@{}}
{#1}--{#2}--\lower3.5ex\vbox{\hbox{{#3}\rule{0ex}{2.5ex}}
\hbox{\hspace{0.4ex}\rule{.2ex}{1ex}\rule{0ex}{1.4ex}}\hbox{{#6}\strut}}--{#4}--{#5}
\end{tabular}}}

\newenvironment{E7}[7]{%
{\small\begin{tabular}{@{}c@{}}
{#1}--{#2}--{#3}--\lower3.5ex\vbox{\hbox{{#4}\rule{0ex}{2.5ex}}
\hbox{\hspace{0.4ex}\rule{.2ex}{1.1ex}\rule{0ex}{1.4ex}}\hbox{{#7}\strut}}--{#5}--{#6}
\end{tabular}}}

\newenvironment{E8}[8]{%
{\small\begin{tabular}{@{}c@{}}
{#1}--{#2}--{#3}--{#4}--\lower3.5ex\vbox{\hbox{{#5}\rule{0ex}{2.5ex}}
\hbox{\hspace{0.4ex}\rule{.2ex}{1ex}\rule{0ex}{1.4ex}}\hbox{{#8}\strut}}--{#6}--{#7}
\end{tabular}}}

\newcommand {\ce}{{\mathfrak c}}

\newcommand {\g}{{\mathfrak g}}
\newcommand {\h}{{\mathfrak h}}

\newcommand {\me}{{\mathfrak m}}

\newcommand {\fN}{\mathfrak{N}}

\newcommand {\q}{{\mathfrak q}}
\newcommand {\rr}{{\mathfrak r}}
\newcommand {\es}{{\mathfrak s}}
\newcommand {\te}{{\mathfrak t}}

\newcommand {\z}{{\mathfrak z}}

\newcommand {\gln}{{\mathfrak{gl}_n}}

\newcommand {\sltn}{{\mathfrak{sl}_{2n}}}
\newcommand {\sltnp}{{\mathfrak{sl}_{2n+1}}}

\newcommand {\spn}{{\mathfrak {sp}}_{2n}}
\newcommand {\son}{\mathfrak{so}_{n}}
\newcommand {\sono}{\mathfrak{so}_{2n+1}}
\newcommand {\sone}{\mathfrak{so}_{2n}}

\newcommand {\gm}{\g\langle m\rangle}
\newcommand {\qm}{\q\langle m\rangle}
\newcommand {\gp}[1]{\g\langle #1\rangle}
\newcommand {\qp}[1]{\q\langle #1\rangle}
\newcommand {\ca}{{\mathcal A}}

\newcommand {\cf}{{\mathcal F}}

\newcommand {\co}{{\mathcal O}}
\newcommand {\cs}{{\mathcal S}}

\newcommand {\BZ}{{\mathbb Z}}
\newcommand {\BN}{{\mathbb N}}

\newcommand {\md}{/\!\!/}

\newcommand{\odin}{{\mathrm{1\hspace{1pt}\!\! I}}}

\newcommand{\lb}{\lambda}
\newcommand{\ap}{\alpha}
\newcommand{\vp}{\varphi}

\renewcommand{\le}{\leqslant}
\renewcommand{\ge}{\geqslant}

\newcommand{\eus}{\EuScript}

\newcommand {\ad}{\mathrm{ad}}
\newcommand {\ads}{{\mathrm{ad}^*}}

\newcommand {\codim}{{\mathrm{codim\,}}}

\newcommand {\gr}{{\mathrm{gr}}}

\newcommand {\ind}{{\mathrm{ind\,}}}
\newcommand {\ins}{\mathrm{ind}}

\newcommand {\rk}{{\mathrm{rk\,}}}

\newcommand {\spe}{\mathrm{Spec}}

\newcommand {\trdeg}{{\mathrm{trdeg\,}}}

\newcommand {\GR}[2]{{\textrm{{\bf #1}}}_{#2}}

\newcommand {\un}{\underline}

\newcommand {\sfb}{\textsl{B\/}}
\newcommand {\sfr}{\mathsf R}
\newcommand {\tinv}{\mathsf{Inv}}
\newcommand {\cgq}{\mathfrak C_\Gamma(\q)}
\newcommand {\ctq}{\mathfrak C_\theta(\q)}

\newfam\eusfam%
\font\Bbbfont=msbm10 scaled 1200%
\def\bbk{\hbox {\Bbbfont\char'174}}

\begin{document}
\setlength{\parskip}{3pt plus 5pt minus 0pt}
\hfill { {\color{blue}\scriptsize \today}}
\vskip1ex

\title[Periodic automorphisms of Takiff algebras]
{Periodic automorphisms of Takiff algebras, contractions, and $\theta$-groups}
\author[D.\,Panyushev]{Dmitri I.~Panyushev}
\address[]{
Independent University of Moscow, Bol'shoi Vlasevskii per. 11,  Moscow 119002, \ Russia
\hfil\break\indent
Institute for Information Transmission Problems, B. Karetnyi per. 19, Moscow 101447
}  
\email{panyush@mccme.ru}
\thanks{Supported in part by  R.F.B.R. grants 05--01--00988
and 06--01--72550.}
\maketitle

\section*{Introduction}
\noindent
Let $G$ be a connected reductive algebraic group with Lie algebra $\g$.
The ground field $\bbk$ is algebraically closed and of characteristic zero.
Fundamental results in invariant theory of the adjoint representation of $G$ 
are primarily associated with 
C.~Chevalley and B.~Kostant. Especially, one should distinguish the 
"Chevalley restriction theorem" and seminal article of Kostant \cite{ko63}. Later, 
Kostant and Rallis  extended these results to the isotropy representation of a symmetric 
variety~\cite{kr71}.
In 1975, E.B.~Vinberg  came up with the theory of $\theta$-groups.
This theory generalises and presents in the most natural form invariant-theoretic 
results previously
known for the adjoint representation and isotropy representations of the symmetric varieties.

Let us remind the main construction and results of Vinberg's article \cite{vi76}.
Let $\theta\in {\rm Aut}(\g)$ be a periodic (=\,finite order) automorphism of $\g$. 
The {\it order\/} of $\theta$ is denoted by $|\theta|$.  
Fix a primitive root of unity $\zeta=\sqrt[|\theta|]1$ and consider 
the {\it periodic grading\/} (or $\BZ_{|\theta|}$-{\it grading})
\[
    \g=\bigoplus_{i\in{\mathbb Z}_{|\theta|}}\g_i \ ,
\]
where $\g_i$ is the $\zeta^i$-eigenspace of $\theta$. In particular, $\g_0=\g^\theta$ is the
fixed point subalgebra for $\theta$.
Let $G_0$ be the connected subgroup of $G$ with Lie algebra  $\g_0$. The restriction
of the adjoint representation yields the natural homomorphism 
$G_0\rightarrow GL(\g_1)$. The linear groups obtained in this way 
are called {\it $\theta$-groups}, and the point is that they have the best possible 
invariant-theoretic properties:

\begin{itemize}
\item \ 
$\bbk[\g_1]^{G_0}$ is a polynomial algebra;
\item \ the quotient morphism  $\pi:\g_1\to \g_1\md G_0=\spe(\bbk[\g_1]^{G_0})$ is flat;
\item \ each fibre of $\pi$ contains finitely many $G_0$-orbits.
\end{itemize}

\noindent
It is a natural problem to extend Vinberg's theory to a more general setting. 
There can be several possible ways for doing so. Here are at least two of them:

(a) \ Determine and investigate a wider class of Lie algebras such that their periodic 
automorphisms lead to representations with similar invariant-theoretic properties. 

(b) \ Given $(\g,\theta)$ as above, construct a non-reductive Lie algebra with good
invariant-theoretic properties, using $\{\g_i\}$ as building blocks.

\noindent These two ways are not mutually exclusive, but in this article we deal with (a).
Our class of Lie algebras consists of Takiff algebras modelled on reductive ones.
Although considerable part of the theory can be developed for arbitrary Takiff algebras, 
substantial applications are related to the reductive case.

For a Lie algebra $\q$ and $m\in\BN$, let $\qm$ denote the Takiff algebra modelled on $\q$.
It is an $\BN$-graded Lie algebra of dimension $m\dim\q$, with nonzero components of
degrees $0,1,\dots,m{-}1$. (See Section~\ref{sect:prelim} for precise definitions.)
Our initial observation is that any periodic $\theta\in{\rm Aut}(\q)$ gives rise to an automorphism
of $\qm$ of the same order, denoted $\widehat\theta$. The fixed point subalgebra of
$\widehat\theta$\/, $\qm^{\widehat\theta}$, \ is a ``mixture'' of eigenspaces of $\theta$, 
i.e., its component of degree $i$ equals $\q_i$, $i=0,1,\dots,m-1$.
Then we consider the
representations of  $\qm^{\widehat\theta}=\qm_0$ 
in eigenspaces of $\widehat\theta$. If $\q$ is quadratic (i.e., $\q\simeq\q^*$), then the coadjoint representation of
$\qm_0$ also occurs in this way.
Our ultimate goal is to describe several instances in which the algebras of invariants
for $(\qm_0,\ad)$ and $(\qm_0,\ads)$
are polynomial. (If the explicit formula for $\q$ is bulky, then we write
$\tinv(\q,\ad)$ (resp. $\tinv(\q,\ads)$) in place of $\bbk[\q]^\q$ (resp. $\cs(\q)^\q$).)
Another observation is that, for special values of $m$, $\qm_0$ is a contraction
of a direct sum of Lie algebras. Namely, 
$\qp{n\vert\theta\vert}_0$ is a contraction of $n\q:=\q\dotplus\ldots\dotplus\q$ 
and $\qp{n\vert\theta\vert+1}_0$ is a contraction of $n\q\dotplus\q_0$.
These are examples of the so-called {\it quasi-graded\/} contractions, and for such 
contractions we establish a rather explicit connection between invariants of two Lie algebras.
For instance, starting from $\bbk[n\q]^{n\q}$, we construct an explicit
subalgebra of $\tinv(\qp{nk}_0,\ad)$, denoted $\eus L_{\bullet}(\bbk[n\q]^{n\q})$.
The graded algebras $\bbk[\q]^\q$ and $\eus L_{\bullet}(\bbk[n\q]^{n\q})$ have the same
Poincar\'e series.
The similar subalgebra for  the coadjoint representation is denoted by 
$\eus L^{\bullet}(\cs(n\q)^{n\q})$.

For $\q=\g$ reductive, we deal with contractions of reductive algebras and therefore 
the theory of $\theta$-groups is also at our disposal. 
Our main result is

\begin{thm}    \label{thm:intro}
Let $\theta$ be a periodic automorphism of  $\g$, $\co^{reg}$ the regular nilpotent 
$G$-orbit in $\g$, and $n\in\BN$  arbitrary. Set $k=\vert\theta\vert$.
\begin{itemize}
\item[\sf (i)] \ Suppose  $\theta$ has the property that $\g_0\cap \co^{reg}\ne\varnothing$. 
Then 
\begin{itemize}
\item \ $\eus L_{\bullet}(\bbk[n\g]^{n\g})=\tinv(\gp{nk}_0, \ad)$ and\/
$\tinv(\gp{nk}_0, \ad)$ is a polynomial algebra of Krull dimension $n{\cdot}\rk\g$.
\item \ $\eus L_{\bullet}(\bbk[n\g\dotplus\g_0]^{n\g\dotplus\g_0})=\tinv(\gp{nk+1}_0, \ad)$ and\/
$\tinv(\gp{nk+1}_0, \ad)$ is a polynomial algebra of Krull dimension $n{\cdot}\rk\g+\rk\g_0$.
\end{itemize}
\item[\sf (ii)] \ Suppose  $\theta$ has the property that $\g_1$ contains regular semisimple 
elements of $\g$ and $\co^{reg}$ meets every irreducible component of the nilpotent 
cone in $\g_1$. Then 
$\eus L^{\bullet}(\cs(n\g)^{n\g})=\tinv(\gp{nk}_0, \ads)$ and\/
 $\tinv(\gp{nk}_0, \ads)$ is a polynomial algebra of Krull dimension $n{\cdot}\rk\g$.
\end{itemize}
\end{thm}
\noindent
The proofs of two parts of this theorem exploit different ideas. Note that 
$\gp{nk+1}_0$ is quadratic and $\gp{nk}_0^*\simeq\gp{nk}_1$. In particular, in part (ii) we describe the invariants for the $\gp{nk}_0$-module $\gp{nk}_1$.

The paper is organised as follows. Sections~\ref{sect:prelim} contains generalities on
quadratic and Takiff algebras and on Lie algebra contractions. 
In Section \ref{sect:quasi}, we introduce quasi-gradings of Lie algebras and 
corresponding contractions. We provide a useful construction of quasi-gradings and  
study the behaviour of invariants.  
In Section~\ref{sect:contr}, we consider periodic automorphisms 
of Takiff algebras and their connection with quasi-graded contractions.
Section~\ref{sect:adj} is devoted to the proof of Theorem~\ref{thm:intro}(i),
and in Section~\ref{sect:coadj} we prove Theorem~\ref{thm:intro}(ii).
Sections~\ref{sect:adj} and \ref{sect:coadj} also contain a number of examples of 
$\theta$-groups that satisfy
the assumptions of Theorem~\ref{thm:intro}. In Section~\ref{sect:end}, we discuss 
open problems and directions for related investigations. 

\subsection*{Notation}  The nilpotent radical of  a Lie algebra $\q$ is denoted by 
$\mathfrak R_u(\q)$. The unipotent radical of an algebraic group $Q$ is $R_u(Q)$.
A direct sum of Lie algebras is denoted with $\dotplus$.  

{\small {\bf Acknowledgements.}  
This work was started during my 
stay at the Max-Planck-Institut f\"ur Mathematik (Bonn) in Spring 2007. 
I am grateful to this institution for warm hospitality and support. }

\section{Preliminaries} 
\label{sect:prelim}

\subsection{Quadratic Lie algebras} \label{sub:quadr}
A Lie algebra $\q$ is called {\it quadratic}, if there is a $\q$-invariant symmetric 
non-degenerate bilinear form on $\q$. Such a form is said to be a {\it scalar product\/} on $\q$.
If $\sfb$ is a scalar product on $\q$, we also say that $(\q,\sfb)$ is quadratic.

Suppose that $(\q,\sfb)$ is quadratic and
$\theta\in{\rm Aut}(\q)$ is  of order $k$. 
Assume that $\sfb$ is a an eigenvector of $\theta$,
i.e., $\sfb(\theta(x),\theta(y))=\zeta^c \sfb(x,y)$ for all $x,y\in\q$ and some $c\in \BZ_k$.
Then $\sfb(\q_i,\q_j)=0$ unless $i+j=c$ (the equality in $\BZ_k$). 
It follows that the dual $\q_0$-module for $\q_i$ is $\q_{c-i}$.
Thus, $\q_0$ is not necessarily  quadratic unless $\zeta^c=1$. 
However,  we  have a weaker
property that the set of $\q_0$-modules $\{\q_i\}$ is closed with respect to taking duals.

If $\q$ is reductive, then $\sfb$ can always be chosen to be $\theta$-invariant, hence $c=0$ 
and $\q_i^*\simeq \q_{-i}$ for all $i\in \BZ_k$. Actually, $\q_0$ is again reductive here.

\subsection{Generalised Takiff  algebras \cite{rt}}   \label{subs:takif}
Let $Q$ be a connected algebraic group with Lie algebra $\q$.
The infinite-dimensional $\bbk$-vector space $\q_\infty:=\q\otimes \bbk[\mathsf T]$ has a 
natural structure of a graded Lie algebra such that 
$[x\otimes \mathsf T^l, y\otimes \mathsf T^k]=[x,y]\otimes \mathsf T^{l+k}$.
Then $\q_{\ge (m)}=\displaystyle\bigoplus_{j\ge m} \q\otimes \mathsf T^j$ is an ideal of
$\q_\infty$, and $\q_\infty/\q_{\ge (m)}$
is a {\it (generalised) Takiff  algebra} (modelled on $\q$), denoted $\q\langle m\rangle$.
Write $Q\langle m\rangle$ for the corresponding connected algebraic group.
Clearly, $\qm$ is $\BN$-graded and $\dim\qm=m\dim\q$. An alternate notation for $\qm$ used
below is
 \[
    \q\ltimes\q\ltimes\ldots\ltimes\q\quad \text{ ($m$ factors)},
 \]
where the consecutive factors from the left to right comes from the subspaces
$\q\otimes T^{i}$, $i=0,1,\dots,m-1$. In particular, 
$\q\langle 2\rangle =\q\ltimes \q$ is the usual semi-direct product.
The image of $\q\otimes T^{i}$ in $\q\ltimes\q\ltimes\ldots\ltimes\q$
 is sometimes denoted by $\q[i]$.
We may (and will) represent  the elements of $\qm$ as vectors:
\[
   \vec x=(x_0,x_1,\dots,x_{m-1}), \text{ where }  \ x_i\in\q[i] .
\]
Recall that $\cs(\q)$ is the symmetric algebra of $\q$ over $\bbk$ and $\bbk[\q]=\cs(\q^*)$.  
For reductive $\g$, $\g[0]$ is a Levi subalgebra of $\gm$ and $\mathfrak R_u({\gm})=
0\ltimes\g\ltimes\dots \ltimes\g$.
It is shown in \cite{rt} that 
$\bbk[\gm]^{\gm}=\bbk[\gm]^{G\langle m\rangle}$ is polynomial.
More precisely, there is an explicit procedure for constructing elements of 
$\bbk[\gm]^{\gm}$ from those of $\bbk[\g]^\g$, which enables us to prove the polynomiality.
The semisimple case with $m=2$ was considered by Takiff in 1971.
In \cite[Sect.\,11]{p05}, we extend results of \cite{rt} to more general Lie algebras and prove
that the algebra $\bbk[\gm]^{R_u(G\langle m\rangle)}$ is also polynomial.

\subsection{Contractions}                            \label{subs:contr}
We only consider Lie algebra contractions of the following form. Let $\mathsf c_t:\q\to \q$, 
$t\in\bbk^\times$, be a polynomial  linear action of $\bbk^\times$ on $\q$. 
That is, $\mathsf c_1={\rm id}$,  $\mathsf c_{t'}\mathsf c_{t''}=\mathsf c_{t'+t''}$,
and the matrix entries of $\mathsf c_{t}$ are polynomials in $t$.
Define the new Lie algebra structure on
the vector space $\q$ by 
\begin{equation}   \label{eq:sdvig-q}
      [x, y]_{(t)}:= \mathsf c_t^{-1}[\mathsf c_t( x), \mathsf c_t( y)] .
\end{equation}
The corresponding Lie algebra is denoted by $\q_{(t)}$. Then $\q_{(1)}=\q$ and
all these algebras are isomorphic. If\/ 
$\lim_{t\to 0}[x, y]_{(t)}$ exists for all $x,y\in\q$, then we obtain a new Lie algebra,
say $\es$, which is a contraction of $\q$. To express this fact, we 
write $\lim_{t\to 0}\q_{(t)}=\es$ or merely $\q\leadsto\es$. We identify $\q$ and $\es$ as vector
spaces.

There is a relation between  invariants of $\q$ and
$\es$. Let $\bbk[\q]_m\simeq \bbk[\q_{(t)}]_m\simeq \bbk[\es]_m$ 
be the space of polynomials of degree $m$ 
and $\bbk[\q_{(t)}]^{\q_{(t)}}_m$ the subspace of invariants of the adjoint
representation of $\q_{(t)}$. 
Then $\dim \bbk[\q_{(t)}]^{\q_{(t)}}_m$ does not depend
on $t\in\bbk^\times$. Therefore $\lim_{t\to 0}\bbk[\q_{(t)}]^{\q_{(t)}}_m$ exists in
the appropriate Grassmannian and  is clearly contained in $\bbk[\es]^\es_m$.
Gathering together components of all degrees, we obtain the subalgebra 
\[     
       \lim_{t\to 0}(\bbk[\q_{(t)}]^{\q_{(t)}})\subset \bbk[\es]^\es .
\]
It can be described more explicitly, as follows. For $F\in \bbk[\q]^{\q}_m$, 
set $F_{(t)}(x)=F(\mathsf c_t(x))$ and 
expand $F_{(t)}=\sum_{j=a}^b F_j{\cdot}t^j$ with $F_a\ne 0$ and $F_b\ne 0$.
We say that $F_a$ (resp. $F_b$) is the {\it initial\/} (resp. {\it highest})
component for $F$ and write $F_\bullet$ (resp. $F^\bullet$) for it. 
Set $\eus L_\bullet(\bbk[\q]^\q_m)=\{F_\bullet \mid F\in \bbk[\q]^\q_m\}$.

\begin{prop}     \label{pr:contr}
Given  $F\in \bbk[\q]^{\q}_m$, we have
\begin{itemize}
\item[{\sf (i)}] \ $F_{(t)}\in \bbk[\q_{(t)}]^{\q_{(t)}}_m$ for any $t\in \bbk^\times$; 
\item[{\sf (ii)}] \  $F_\bullet \in (\lim_{t\to 0} \bbk[\q_{(t)}]^{\q_{(t)}}_m)
\subset \bbk[\es]^\es_m$;
\item[\sf (iii)] \ $\lim_{t\to 0} \bbk[\q_{(t)}]^{\q_{(t)}}_m
=\eus L_\bullet(\bbk[\q]^\q_m)$.
\end{itemize}
\end{prop}
\begin{proof} (i) Consider the representation of $\q$ in $\bbk[\q]$.  
Let $e\ast (?)$ denote
the operator corresponding to $e\in\q$. Then $F\in\bbk[\q]^\q$ if and only if 
$e\ast F=0$ for any $e\in\q$. The corresponding operator for $e\in\q_{(t)}$ 
is denoted by $e\ast_{(t)}(? )$.
Note that $\mathsf c_t:\q\to \q$ can be regarded as isomorphism
of Lie algebras $\q$ and $\q_{(t)}$.
Therefore  $(e\ast_{(t)} F_{(t)})(x)=(\mathsf c_t(e)\ast F)(\mathsf c_t(x))$,
which yields the assertion.

(ii) \  Consider $F_{(t)}$ as a curve in the projectivisation of $\bbk[\q]_m$.
Then $\lim_{t\to 0}F_{(t)}=F_\bullet$.

(iii) \ The inclusion "$\supset $" is already proved.
 Hence it suffices to show that
$\dim \bbk[\q]^{\q}_m=\dim \eus L_\bullet(\bbk[\q]^\q_m)$.
Consider the finite ascending filtration  $\cf_\bullet$ of $\bbk[\q]^\q_m$:
\[
   \bbk[\q]^\q_m=\cf_0\supset \dots \supset \cf_j\supset \dots ,
\]
where
\[
     \cf_j:=\{F\in \bbk[\q]^\q_m  \mid  F_{(t)}=\sum_{j\ge i} F_j t^j \} .
\]
Let $\mathcal F_N$ be the last nonzero term of $\cf_\bullet$.
Set ${\rm Gr}(\bbk[\q]^\q_m)=\bigoplus_{i=0}^N \cf_i/\cf_{i+1}$.
Then $\dim \bbk[\q]^\q_m=\dim {\rm Gr}(\bbk[\q]^\q_m)$, and for each $i$ there  is a 
natural linear mapping
$\vp_i: \cf_i/\cf_{i+1} \to \eus L_\bullet(\bbk[\q]^\q_m)$,
$F+\cf_{i+1} \mapsto F_\bullet$.
 Clearly, each $\vp_i$ is injective.
Furthermore, $\vp_i(\cf_i/\cf_{i+1})$ consists of polynomials of weight $(-i)$ with respect
to the induced action of $\bbk^\times$ in $\bbk[\q]_m$. Hence, the subspaces 
$\vp_i(\cf_i/\cf_{i+1})$ are linearly independent.
\end{proof}
\noindent 
Similar results hold for coadjoint representations and algebras of invariants 
$\cs(\q_{(t)})^{\q_{(t)}}$.
The only notable difference is that for $F\in \cs(\q)^\q$ one have to take the highest component
$F^\bullet$.

\section{Quasi-graded contractions and invariants} 
\label{sect:quasi}

\noindent 
A Lie algebra $(\q, [\ ,\ ])$ is said to be {\it quasi-graded\/} if there is a vector space decomposition
$\q=\bigoplus_{i=0}^{k-1}\q_i$ such that  $[\q_i,\q_j]\subset \q_{i+j}$ whenever $i+j\le k-1$.
There are no conditions on $[\q_i,\q_j]$ if $i+j\ge k$. The family of
subspaces $\Gamma=\{\q_i\}_{i=0}^{k-1}$ is said to form a 
{\it quasi-graded structure (of order $k$)} on $\q$. 
Define the new Lie algebra structure,
$[\ ,\ ]_{\Gamma}$, on the vector space $\q$ as follows:
for $x_i\in\q_i$, we set \ 
$[x_i,x_j]_{\Gamma}:=\begin{cases} [x_i,x_j],  &\text{if }\ i{+}j\le k{-}1; \\
0,  & \text{if }\ i{+}j \ge k. \end{cases}$. \ The resulting $\BN$-graded Lie algebra is 
denoted by 
$\cgq$ or $\q_0\ltimes\q_1\ltimes\ldots\ltimes\q_{k-1}$.  
It is a contraction of $\q$ in the sense of Subsection~\ref{subs:contr}. 
Indeed, for $t\in\bbk^\times$, 
define  linear operators $\mathsf c_t:\q\to \q$  by
$\mathsf c_t\vert_{\q_i}=t^i{\cdot}{\rm id}$ and define $\q_{(t)}$ as above.
Then  
\[
[x_i,x_j]_{(t)}=\begin{cases} [x_i,x_j],  &\text{if }\ i{+}j\le k{-}1; \\
t{\cdot}(\text{\rm a polynomial expression in $t$}),  & \text{if }\ i{+}j \ge k. \end{cases} 
\] 
It follows that $\lim_{t\to 0} \q_{(t)}
=\cgq$. The passage $\q\leadsto \cgq$ as well as $\cgq$ itself
is called the $\Gamma$-{\it contraction\/}  or a {\it quasi-graded contraction\/} of $\q$. 
We identify $\q$ and $\cgq$ as graded vector spaces. 

\begin{ex}  Here we provide important examples of quasi-graded structures.

1)  Let $\theta\in{\rm Aut}(\q)$ be of order $k$. Then the corresponding
$\BZ_k$-grading of $\q$ is a quasi-grading. The quasi-graded contraction 
associated with $\theta$ is said to be 
{\it cyclic\/} or a $\BZ_k$-{\it contraction}. It is denoted by
$\ctq$.

2) Let $\h$ be a subalgebra of $\q$. Suppose that there is an
$\ad(\h)$-stable subspace $\me\subset \q$ such that $\q=\h\oplus\me$. 
Then $\Gamma=\{\h,\me\}$ is a quasi-graded structure of order 2 on $\q$. The passage 
$\q\leadsto \h\ltimes\me$ is called the $\h$-{\it isotropy contraction\/} of $\q$.
Here $\h\ltimes\me$ is  the {\it semi-direct product\/} of
(the Lie algebra) $\h$  and (the $\h$-module) $\me$.
\end{ex}

\begin{rem}
Any $\BZ_2$-contraction is an isotropy contraction;
$\BZ_2$-contractions of reductive Lie algebras have already been considered in
\cite[\S\,9]{p05} and \cite{coadj}.
\end{rem}

There is a natural construction of certain quasi-gradings from periodic gradings.

\begin{lm}     \label{lm:quasi}
Suppose $\theta\in{\rm Aut}(\q)$ is of order $k$. Let\/ $\q_0$ be the corresponding fixed point subalgebra of $\q$. Then
$\q_0\dotplus\q$ has a natural quasi-graded structure of order $k+1$.
The corresponding quasi-graded contraction is isomorphic to
$\q_0\ltimes\q_1\ltimes\ldots\ltimes\q_{k-1} \ltimes\q_0$.
\end{lm}\begin{proof}
Consider the family of subspaces $\Gamma= \{\es_j\}_{j=0}^k$ such that
$\es_0=\{(x,x)\mid x\in\q_0\}\subset \q_0\dotplus\q_0\subset \q_0\dotplus \q$, 
$\es_i=\q_i\subset\q$ for $i=1,\dots,k-1$, and $\es_k=\q_0\subset\q$.
Then $\q_0\dotplus\q=\bigoplus_{i=0}^k\es_i$ and
\begin{equation}   \label{eq:non-cyclic}
[\es_i,\es_j]\subset \begin{cases}  \es_{i+j}, & \text{if }\ i+j\le k, \\
 \es_{i+j-k}, & \text{if }\ i+j\ge k+1. \end{cases}
 \end{equation} 
Therefore $\Gamma$ is a quasi-grading. Obviously,
$\mathfrak C_{\Gamma}(\q_0\dotplus\q)\simeq \q_0\ltimes\q_1\ltimes\ldots\ltimes\q_{k-1}
\ltimes\q_0$.
\end{proof}

\begin{rem}   The quasi-graded contraction in the lemma is not cyclic.
\end{rem}

\noindent
Below, we explicitly describe a connection between invariants of $\q$ and 
$\mathfrak C_{\Gamma}(\q)$. In the special case of $\BZ_2$-contractions of reductive Lie
algebras, such a connection is explained in \cite[Proposition\,3.1]{coadj}.
Discussing ``invariants of $\q$", we always mean  invariants of the
adjoint and coadjoint representations, i.e., the algebras $\bbk[\q]^\q$ and
$\cs(\q)^\q$. If an explicit formula for $\q$ appears to be bulky, we also use notation
$\tinv(\q, \ad)$ and $\tinv(\q, \ads)$.

Any quasi-grading $\Gamma=\{\q_i\}_{i=0}^{k-1}$ determines a polygrading of $\cs(\q)$.
Using the vector space
decomposition $\q=\bigoplus_{j=0}^{k-1} \q_j$, any $F\in\cs(\q)$ can be written as
$F=\sum F_{i_0,i_1,\dots,i_{k-1}}$, where 
$F_{i_0,i_1,\dots,i_{k-1}}\in \bigotimes_{j=0}^{k-1}\cs^{i_j}(\q_j)$.
Define a specialisation of this polygrading by
\[
     \deg_\Gamma(F_{i_0,i_1,\dots,i_{k-1}})=i_1+2i_2+\dots +(k-1)i_{k-1}.
\]
We say that $\deg_\Gamma$ is the $\Gamma$-{\it degree} in $\cs(\q)$.
Note that the usual degree in $\cs(\q)$ is defined by 
$\deg(F_{i_0,i_1,\dots,i_{k-1}})=i_0+i_1+\dots +i_{k-1}$.
If we refer to homogeneous polynomials,  then the usual degree is meant,
unless otherwise stated. 
If $F$ is homogeneous, then $\deg(F)=\deg(\gr^\bullet F)=\deg(\gr_\bullet F)$.
Using the dual decomposition of $\q^*$, one similarly defines the $\Gamma$-degree for
polynomials in $\bbk[\q]$.

\begin{df} 
For $F\in\cs(\q)$, let
$\gr^\bullet F$ (resp. $\gr_\bullet F$) denote the component of $F$ of the {\it maximal\/} 
(resp. {\it minimal}) $\Gamma$-degree. The same notation applies to $F\in\bbk[\q]$.
\end{df}

\noindent
As $\q$ and $\cgq$ are identified as vector spaces,
each $F\in\cs(\q)$ (resp. $F\in\bbk[\q]$) can also be regarded as element of $\cs(\cgq)$
(resp. $\bbk[\cgq]$).

\begin{thm}    \label{thm:inv-contr}
Let $\Gamma=\{\q_i\}_{i=0}^{k-1}$ be an arbitrary quasi-graded structure on $\q$. 
\begin{itemize}   
\item[\sf (i)] \  If  $F\in \cs(\q)^{\q}$, then $\gr^\bullet F \in \tinv(\cgq,\ads)$.
\item[\sf (ii)] \  
If  $F\in \bbk[\q]^{\q}$, then $\gr_\bullet F \in \tinv(\cgq,\ad)$.
\end{itemize}
\end{thm}\begin{proof}
(i) 
Let $\{\ ,\ \}$ (resp. $\{\ ,\ \}_{\Gamma}$) denote the Poisson bracket in $\cs(\q)$
(resp. $\cs(\cgq)$). 
As is well-known, $F\in\cs(\q)^\q$ if and only if $\{x,F\}=0$ for any $x\in\q$.

If $x_i\in\q_i$, then $\deg_\Gamma(x_i)=i$.  
Therefore 
\[
\deg_\Gamma\{x_i,x_j\}\,\begin{cases} =i+j,& \text{if }\ i+j< k,\\
<i+j, & \text{if }\ i+j\ge k.\end{cases}
\]
If $\deg_\Gamma(F_{i_0,i_1,\dots,i_{k-1}})=m$, then 
$\deg_\Gamma(\{x_j, F_{i_0,i_1,\dots,i_{k-1}}\}_{\Gamma})=m+j$. Furthermore,
comparing the commutators in $\q$ and $\cgq=
\q_0\ltimes\ldots\ltimes\q_{k-1}$ shows that 
\[
    \{x_j, F_{i_0,i_1,\dots,i_{k-1}}\}=\{x_j, F_{i_0,i_1,\dots,i_{k-1}}\}_{\Gamma}+ 
    (\text{terms of  $\Gamma$-degree $< j+m$}).
\]  
It follows that  
$\{x_j,\gr^\bullet F\}_{\Gamma}$ is the component of the maximal possible 
$\Gamma$-degree in $\{x_j,F\}$.
As $\{x_j,F\}=0$, we also must have $\{x_j,\gr^\bullet F\}_{\Gamma}=0$.

(ii) This follows from Proposition~\ref{pr:contr}, since
$\gr_\bullet F$ is the initial component of $F_{(t)}$ with respect to $t$.
\end{proof}

\begin{rem}
Proposition~\ref{pr:contr} can be adapted for obtaining another proof of part (i). In the context
of quasi-graded contractions, we prefer to use notation $\gr_\bullet F$ in place of $F_\bullet$
(and likewise for $\gr^\bullet F$). 
\end{rem}

\noindent
Recall that $\eus L_\bullet(\bbk[\q]^\q):=\{\gr_\bullet F\mid F\in \bbk[\q]^\q\}$ is a subalgebra
of $\bbk[\cgq]^{\cgq}$, and it follows from Proposition~\ref{pr:contr}(iii) that the graded 
algebras $\bbk[\q]^\q$ and $\eus L_\bullet(\bbk[\q]^\q)$ have equal Poincar\'e series.
Similarly, $\eus L^\bullet(\cs(\q)^\q)$ is a subalgebra of $\cs(\cgq)^{\cgq}$.

\section{Periodic automorphisms of Takiff algebras and related contractions } 
\label{sect:contr}

\noindent
Throughout this section, $\theta$ is a periodic automorphism of $\q$ and
$\zeta=\sqrt[|\theta|]1$; usually $|\theta|=k$.

\noindent
Every periodic  $\theta\in{\rm Aut}(\q)$ can be extended to an automorphism
of $\qm$ of the \un{same} order.
There are at least two ways for doing so:
\begin{gather}
\notag \widehat\theta\vert_{\q[i]}=\theta ;  \\  \label{eq:case-b}
\widehat\theta\vert_{\q[i]}=\zeta^{-i}\theta .
\end{gather}
In both cases, it is easily seen that  $\widehat\theta\in {\rm Aut}(\q\langle m\rangle)$. 
Notice that no relation between $m$ and $|\theta|$ is required! In the first case, the fixed point
subalgebra is again a Takiff algebra. Therefore we only work with the second case,
which is certainly more interesting. 
That is, from now on,  
$\widehat\theta$ is defined by Eq.~\eqref{eq:case-b}.
Then the $\zeta^i$-eigenspace of $\widehat\theta$\/ is 
\[
\q\langle m\rangle_i=
  \q_i\ltimes\q_{i+1}\ltimes\ldots\ltimes\q_{i+m-1} ,
\] 
where $\q_i$ is the $\zeta^i$-eigenspace of $\theta$ and
all subscripts are regarded as elements of $\BZ_{\vert\theta\vert}$. 
In particular, the fixed point subalgebra for $\widehat\theta$ is
\[  
\qm^{\widehat\theta}=\q\langle m\rangle_0=
  \q_0\ltimes\q_1\ltimes\ldots\ltimes\q_{m-1} . 
\]
If $(\q,\sfb)$ is quadratic, then $\qm$ is also a quadratic Lie algebra.
For,  we can extend $\sfb$ to $\qm$ by letting
\begin{equation}    \label{eq:quadr}
    \hat\sfb(\vec x,\vec y)=\sum_{i=0}^{m-1} \sfb(x_i, y_{m-1-i}) \ .
\end{equation}

\begin{lm}   \label{lm:value}
Suppose $(\q,\sfb)$ is quadratic and $\sfb$ is a $\theta$-eigenvector with eigenvaue
$\zeta^c$. Then $\hat\sfb$ is  a
$\widehat\theta$-eigenvector with eigenvalue  $\zeta^{c+1-m}$.
\end{lm}\begin{proof} This follows from Eq.~\eqref{eq:case-b} and \eqref{eq:quadr}.
\end{proof}

\noindent Thus,  even if a scalar product on $\q$ is $\theta$-invariant,
then extending it to a Takiff algebra $\qm$ we may obtain the scalar product 
with a non-trivial $\widehat\theta$-eigenvalue.
Combining Lemma~\ref{lm:value} and Subsection~\ref{sub:quadr}, we obtain

\begin{cl}    \label{cl:c}
If $(\q,\sfb)$ is quadratic and $\sfb$ is a $\theta$-eigenvector with eigenvaue
$\zeta^c$, then the dual $\qm_0$-module for $\qm_i$ is $\qm_{c+1-m-i}$. 
In particular, if $c=0$, then both
$(\q_0, \sfb\vert_{\q_0})$  and $(\qp{nk+1}_0, \hat\sfb)$, $n\in\BN$, are also quadratic.
\end{cl}

\noindent
We are going to study the algebra of invariants for the adjoint and coadjoint 
representations of algebras of the form $\qm_0$. Since
$\qm_0^*\simeq \qm_{c+1-m}$, 
the coadjoint representation of $\qm_0$ can also be realised as 
$\widehat\theta$-representation.

\noindent
Let us point out  two important cases:
\begin{gather*}  
\text{If $m=k$, then $\qp{k}_0=\q_0\ltimes\q_1\ltimes\ldots\ltimes\q_{k-1}$,}    
\\ 
\text{If $m=k+1$, then $\qp{k+1}_0= \q_0\ltimes\q_1\ltimes\ldots\ltimes\q_{k-1}
\ltimes\q_0$.}
\end{gather*} 
In the first case, each $\q_i$ occurs exactly once.
Furthermore, if $c=0$, then $\qp{k}_0^*\simeq \qp{k}_1$ and
$\gp{k+1}_0$  is quadratic.
The utility of these (and some other related) cases is that the fixed point subalgebra for
$\widehat\theta$ appears to be a quasi-graded contraction.

\begin{prop}   \label{prop:contr}
{\sf (i)} \ The passage $\q\leadsto \qp{k}_0$ is the $\BZ_k$-contraction associated with
$\theta$;  {\sf (ii)} \
The passage\/ $\q_0\dotplus\q\leadsto \qp{k+1}^{\widehat\theta}=
\qp{k+1}_0$ is a quasi-graded contraction of order\/
$k+1$.
\end{prop}\begin{proof}
(i) is obvious; (ii) follows from Lemma~\ref{lm:quasi}.
\end{proof}

\noindent
Recall that the cyclic contraction of $\q$ associated with $\theta$ is denoted by 
$\mathfrak C_\theta(\q)$.
Therefore Proposition~\ref{prop:contr}(i) can be expressed as the equality
$\mathfrak C_\theta(\q)\simeq \qp{|\theta|}^{\widehat\theta}$.
To generalise previous observations, we need some preparations. 

\begin{lm}    \label{lm:ext-theta}
Let $\theta$ be a periodic automorphism of $\q$. For any $n\in\BN$, there
is a periodic automorphism $\widetilde\theta$ of\/ $n\q:=\q\dotplus\ldots\dotplus\q$ 
($n$ summands) such that  $|\widetilde\theta|=n|\theta|$.
Furthermore, for any $j\in \BZ_{n|\theta|}$, we have  $(n\q)_j\simeq
\q_{\bar j}$, where $\bar j$ is the image of $j$ in $\BZ_{|\theta|}$.
\end{lm}\begin{proof}
Let $k=|\theta|$ and $\zeta=\sqrt[k]1$. Using the  $\BZ_k$-grading
$\q=\bigoplus_{i\in\BZ_k}\q_i$, we consider $n\q$ as the direct
sum of spaces $n(\q_i)=\q_i\dotplus\ldots\dotplus\q_i$. All these spaces are to be
$\widetilde\theta$-stable,  and we define $\widetilde\theta$ for each $i$ separately.
For $(a_1,\dots,a_n)\in n(\q_i)$,  set
\[
     \widetilde\theta(a_1,\dots,a_n)=(a_2,\dots,a_n,\zeta^ia_1) .
\]
Obviously, $\widetilde\theta\in {\rm Aut}(n\q)$  and $\widetilde\theta^n\vert_{n(\q_i)}=
\zeta^i{\cdot}\text{id}$. Hence $|\widetilde\theta|=nk$.
To describe the eigenspaces of $\widetilde\theta$, we choose a primitive
root $\mu=\sqrt[nk]1$ such that $\mu^n=\zeta$. Write $(n\q)_j$ for the $\mu^j$-eigenspace of
$\widetilde\theta$. We claim that
\[
  (n\q)_j=\{x,\mu^j x,\dots,\mu^{(n-1)j}x \mid x\in \q_{\bar j}\}.
\]
Indeed, let $\rr_j$ denote the right-hand side.
It is easily seen that $\rr_j\subset (n\q)_j$,  
$\rr_j\simeq \q_{\bar j}$, and $\rr_j\subset n(\q_{\bar j})$. 
Moreover,  for any $j\in\{0,1,\dots,k-1\}$, the sum $\sum_{l=0}^{n-1}\rr_{j+lk}$ is direct
(use the Vandermonde determinant!).
Whence $n\q=\bigoplus_{j=0}^{nk-1} \rr_j$, and we are done.
\end{proof}
Our general result on contractions associated with Takiff algebras is

\begin{thm}   \label{thm:main-contr} 
Given a periodic $\theta\in{\rm Aut}(\q)$ and  $n\in\BN$, 
consider $\widetilde\theta\in{\rm Aut}(n\q)$ as above.  
Then
\begin{itemize}
\item[\sf (i)] \ For $\widehat\theta\in{\rm Aut}(\qp{n|\theta|}$, the fixed point subalgebra
$\qp{n|\theta|}^{\widehat\theta}$ is the cyclic contraction of $n\q$ associated with 
$\widetilde\theta$, i.e., $\qp{n|\theta|}^{\widehat\theta} \simeq 
\mathfrak C_{\widetilde\theta}(n\q)$;
\item[\sf (ii)] \  For $\widehat\theta\in{\rm Aut}(\qp{n|\theta|+1}$, 
the fixed point subalgebra $\qp{n|\theta|+1}^{\widehat\theta}$ is a 
quasi-graded contraction (of order $nk+1$) of  $\q_0\dotplus n\q$.
\end{itemize}
\end{thm}\begin{proof}
This is an immediate consequence of  our previous results.  First,
we consider the $\BZ_{nk}$-grading of $n\q$ constructed in Lemma~\ref{lm:ext-theta},
which yields (i). For (ii), we endow $\q_0\dotplus n\q$ with the quasi-graded structure 
of order $nk+1$ using Lemmas~\ref{lm:quasi} and  \ref{lm:ext-theta}. We also need 
the fact that $(n\q)_j\simeq \q_{\bar j}$ for any $j=0,1,\dots,nk-1$.
\end{proof}

\begin{cl}      \label{cl:cyc_tak}
Let $\widetilde\theta$ be the cyclic permutation of the summands in $n\q$. Then
$ \qp{n} \simeq \mathfrak C_{\widetilde\theta}(n\q)$.
\end{cl}\begin{proof}
This is the particular case of Theorem~\ref{thm:main-contr}(i) with $\theta={\rm id}$ and 
hence $\widehat\theta={\rm id}$.
\end{proof}

\noindent
Note that Proposition~\ref{prop:contr} corresponds to the case $n=1$ in 
Theorem~\ref{thm:main-contr}.

\begin{rmk}  \label{rem:future}
The $\BZ_{nk}$-contraction of $n\q$ is the $\BN$-graded algebra
$(n\q)_0\ltimes(n\q)_1\ltimes\ldots \ltimes(n\q)_{nk-1}$, but using the isomorphisms 
$(n\q)_j\simeq \q_{\bar j}$ we can write it as
$\qp{nk}^{\widehat\theta}=\q_0\ltimes\q_1\ltimes\ldots \ltimes\q_{k-1}$, 
where each $\q_i$ occurs $n$ times.
Yet, one should not forget that different summands $\q_i$ in the last expression
arise from different subspaces of $n\q$.
For future use, we record the fact that the subalgebra $\q_0$ in 
$\qp{nk}^{\widehat\theta}$ or $\qp{nk+1}^{\widehat\theta}$
(the very first summand) corresponds under the isomorphisms of Lemma~\ref{lm:ext-theta}
and Theorem~\ref{thm:main-contr} to the diagonally embedded subalgebra $\q_0$ in 
$n(\q_0)\subset n\q$ or in $(n+1)\q_0\subset \q_0\dotplus n\q$.
\end{rmk}

Combining Theorems~\ref{thm:inv-contr} and \ref{thm:main-contr}, we obtain

\begin{thm}    \label{thm:inv-contr2}     
Let $\theta\in{\rm Aut}(\q)$ be periodic and $n\in\BN$. Consider the graded structure of 
$n\q$ and quasi-graded structure of $\q_0\dotplus n\q$ determined by $\theta$.
\begin{itemize}   
\item[\sf (i)] \ If $F\in \tinv(n\q, \ads)$, 
then $\gr^\bullet F \in \tinv(\qp{n|\theta|}_0,\ads)$;
\item[\sf (ii)] \ If $F\in \tinv(n\q, \ad)$, 
then $\gr_\bullet F \in \tinv(\qp{n|\theta|}_0,\ad)$.
\item[\sf (iii)] \ If $F\in \tinv(\q_0\dotplus n\q, \ads)$, 
then $\gr^\bullet F \in \tinv(\qp{n|\theta|+1}_0,\ads)$;
\item[\sf (iv)] \ If $F\in \tinv(\q_0\dotplus n\q, \ad)$, 
then $\gr_\bullet F \in \tinv(\qp{n|\theta|+1}_0,\ad)$.
\end{itemize}
\end{thm}

\noindent
Let us point out  possible applications of these procedures.  For definiteness, consider
part (ii).
Taking each (homogeneous) $F\in \bbk[n\q]^{n\q}\simeq
\tinv (\q,\ad)^{\otimes n}$ to $\gr_\bullet F$, 
we obtain the subalgebra of $\tinv(\qp{n|\theta|}_0,\ad)$, which is denoted by 
$\eus L_\bullet (\bbk[n\q]^{n\q})$.
We know that $\bbk[n\q]^{n\q}$ and $\eus L_\bullet (\bbk[n\q]^{n\q})$ 
have the same Poincar\'e series with respect to the usual degree. 
If one knows somehow that 
$\bbk[n\q]^{n\q}$ and $\tinv(\qp{n|\theta|}_0,\ad)$ also have equal Poincar\'e series, then
we obtain the equality $\eus L_\bullet (\bbk[n\q]^{n\q})=\tinv(\qp{n|\theta|}_0,\ad)$.
This sometimes allows us to prove that good properties of $\tinv (\q,\ad)$
are carried to $\tinv (\qp{nk}_0,\ad)$ over.
Instances of such a phenomenon are discussed in the following sections.

\section{Good cases for adjoint representations} 
\label{sect:adj}

\noindent
From now on, our initial object is a reductive Lie algebra $\g$
(in place of an arbitrary $\q$). As above, it is assumed that 
$\theta\in{\rm Aut}(\g)$  and $|\theta|=k$. 
If an algebra of invariants appears to be  graded polynomial, then the 
elements of any set of algebraically independent homogeneous generators 
will be referred to as {\it basic invariants}. E.g., one can consider basic invariants in 
$\cs(\g)^\g=\cs(\g)^G$ or $\bbk[\g_1]^{G_0}$.

Let us start with some remarks concerning the reductive case.
We assume that the scalar product on $\g$ is $\theta$-invariant. Hence $c=0$ in 
Corollary~\ref{cl:c} and
$\gp{nk+1}_0$ is quadratic for any $n\in\BN$.
Recall that $\g[0]\simeq \g$ is a  Levi subalgebra of $\gm$ and 
$\mathfrak R_u(\gm)=0\ltimes \g\ltimes\ldots\ltimes\g$. 
Similarly, $\g_0\subset \g[0]$ is a Levi subalgebra of $\gm_0$ and
$\mathfrak R_u(\gm_0)=0\ltimes \g_1\ltimes\ldots\ltimes\g_{m-1}$.

Let $\fN$ be the nilpotent cone of $\g$ and 
$\co^{reg}\subset\fN$  the regular nilpotent $G$-orbit.

\begin{thm}   \label{teor1+2}
Suppose  $\theta$ has the property that $\g_0\cap \co^{reg}\ne\varnothing$ 
and $n\in\BN$ is arbitrary. Then 
\begin{itemize}
\item[\sf (i)] \ 
$\eus L_{\bullet}(\bbk[n\g]^{n\g})=\tinv(\gp{nk}_0, \ad)$ and\/
$\tinv(\gp{nk}_0, \ad)$ is a polynomial algebra of Krull dimension $n{\cdot}\rk\g$.
\item[\sf (ii)] \ 
$\eus L_{\bullet}(\bbk[n\g\dotplus\g_0]^{n\g\dotplus\g_0})=\tinv(\gp{nk+1}_0, \ad)$
and $\tinv(\gp{nk+1}_0, \ad)$ is a polynomial algebra of Krull dimension 
$n{\cdot}\rk\g+\rk\g_0$.
\end{itemize}
\end{thm}\begin{proof}
(i) \ Consider the chain of Lie algebra contractions:
\[
  n\g \underset{(1)}{\leadsto}  \gp{nk}_0=
  \g_0\ltimes\g_1\ltimes\ldots\ltimes\g_{k-1} 
  \underset{(2)}{\leadsto} 
   \g_0\ltimes(\underbrace{\g_1\oplus\ldots \oplus\g_{k-1}}_{nk-1})=:\h
   \ ,
\]
where ``$\underset{(1)}{\leadsto}$'' is constructed in Theorem~\ref{thm:main-contr}(i) and 
``$\underset{(2)}{\leadsto}$'' stands for the $\g_0$-isotropy contraction of $\gp{nk}_0$. 
In other words, $\mathfrak R_u(\gp{nk}_0)$ is proclaimed to be abelian in $\h$.
The direct passage $n\g\leadsto \h$ can also be regarded as the $\g_0$-isotropy 
contraction of $n\g$, where $\g_0$ is the diagonally embedded subalgebra 
$\g_0\subset n(\g_0)\subset n\g$
(cf. Remark~\ref{rem:future}).
Consider the corresponding transformations of algebras of invariants
\begin{equation}   \label{eq:transform}
   \bbk[n\g]^{n\g}\leadsto \eus L_\bullet^{(1)}(\bbk[n\g]^{n\g})\hookrightarrow \tinv(\gp{nk}_0,\ad)\leadsto
   \eus L_\bullet^{(2)}(\tinv(\gp{nk}_0,\ad))\hookrightarrow \bbk[\h]^\h .
\end{equation}
Here the functor $\eus L_\bullet^{(i)}$ corresponds to the contraction 
$\underset{(i)}{\leadsto}$ and each 
arrow "$\leadsto$" preserves the Poincar\'e series. 
By \cite[Theorem\,6.2]{p05}, the algebra $\bbk[\h]^\h$ is  polynomial. That is, 
both extreme algebras are polynomial. 
Furthermore, the degrees of the basic invariants  for both
are the same \cite[Theorem\,9.5(2)]{p05}. (The crucial property is that the Levi subalgebra
$\g_0\subset \h$ arises from the diagonally embedded $\g_0\subset n(\g_0)\subset n\g$,
and therefore a regular nilpotent element of  $\g_0$ is still regular in $n\g$.)

The equality for the degrees implies that both embeddings in
Eq.~\eqref{eq:transform}  are actually isomorphisms, and all the algebras involved have
one and the same Poincar\'e series. 
Let $F_1,\dots,F_l$ be the basic invariants in
$\bbk[\h]^{\h}$. Then there are homogeneous $\tilde F_i\in \bbk[n\g]^{n\g}$ such that 
$\gr_\bullet^{(2)} (\gr_\bullet^{(1)}(\tilde F_i))=F_i$ for all $i$. 
It is then easily seen that $\{\tilde F_i\}$ and  $\{\gr_\bullet^{(1)}\tilde F_i\}$ must 
be algebraically independent. Since $\deg F_i=\deg(\gr_\bullet^{(1)}\tilde F_i)$, 
$\tinv(\gp{nk}_0,\ad)$ is freely generated by 
the $\{\gr_\bullet^{(1)}\tilde F_i\}$'s.


(ii) We begin with the chain of Lie algebra contractions:
\[
  \g_0\dotplus n\g \underset{(1)}{\leadsto}  \gp{nk+1}_0=
  \g_0\ltimes\g_1\ltimes\ldots\ltimes\g_{k-1}\ltimes\g_{0} 
  \underset{(2)}{\leadsto} 
   \g_0\ltimes(\underbrace{\g_1\oplus\ldots \oplus\g_{k-1}\oplus\g_{0}}_{nk})
   \ .
\]
The rest is essentially the same as in part (i). 
\end{proof}

The list of $\theta$ satisfying the assumption $\g_0\cap\co^{reg}\ne\varnothing$ is not long. 
For simple $\g$, all possible pairs $\g\supset\g_0$ are pointed out below:
\begin{equation}   \label{list}
\begin{split}
|\theta|=2: & \quad \sltn\supset\spn,\ \sltnp\supset\sono,\ \sone\supset \mathfrak{so}_{2n-1},\  
  \GR{E}{6}\supset \GR{F}{4};\\
|\theta|=3: & \quad \GR{D}{4}\supset \GR{G}{2}.
\end{split}
\end{equation}

\noindent
For semisimple Lie algebras, the only new possibility is the cyclic permutation of summands in
$n\g$, which leads to $\gp{n}$, cf. Corollary~\ref{cl:cyc_tak}.

\begin{ex}   \label{real-life}
We give  realisations of algebras 
$\gp{nk}_0$ and $\gp{nk+1}_0$, associated with the list in~\eqref{list}, as centralisers of 
nilpotent elements. If $\theta$ is an involution, then there are only two eigenspaces, $\g_0$ and $\g_1$, and
we will use  the more suggestive notation $\mathfrak L_m(\g_0,\g_1)$ in place of $\gm_0$. 
For $\g$ simple, $\mathfrak L_m(\g_0,\g_1)$ is quadratic if and only if $m$ is odd.
Irreducible $\g_0$-modules occurring in $\g_1$ are depicted by their
highest weights. Namely, $\sfr(\lb)$ is a simple module with highest weight $\lb$.
The $i$-th fundamental weight of a simple Lie algebra is denoted by $\varpi_i$, with 
the numbering from \cite{vo}. We write $\odin$ for the trivial 1-dimensional module.
The symbol $\z_n$ stands for the $n$-dimensional centre.
Items 1$^o$--4$^o$ below provide realisations of the algebras 
$\mathfrak L_m(\g_0,\g_1)$ associated with the {\it outer\/} involutions of 
$\g=\mathfrak{sl}_N$ or $\mathfrak{gl}_N$.
\begin{itemize}
\item[1$^o$.]  Let  $e\in \tilde\g=\mathfrak{sp}_{2nm}$ be a nilpotent element 
with partition $(\underbrace{2m,2m,\dots,2m}_{n})=:((2m)^{n})$. 
Then $\tilde\g^e\simeq \mathfrak L_{2m}(\son,\sfr(2\varpi_1)\oplus \odin)\simeq 
\mathfrak L_{2m}(\son,\sfr(2\varpi_1))\dotplus \z_m$.
Here $\son\oplus(\sfr(2\varpi_1)\oplus \odin)=\gln$ and $\tilde\g^e$ is a contraction of
$m(\gln)$.
\item[2$^o$.]  Let  $e\in \tilde\g=\mathfrak{so}_{n(2m+1)}$ be a nilpotent element 
with partition $((2m+1)^{n})$. 
Then $\tilde\g^e\simeq\mathfrak L_{2m+1}(\son,\sfr(2\varpi_1)\oplus \odin)\simeq 
\mathfrak L_{2m+1}(\son,\sfr(2\varpi_1))\dotplus \z_m$.
Here  $m(\gln)\dotplus\son\leadsto \tilde\g^e$.
\item[3$^o$.]  Let  $e\in \tilde\g=\mathfrak{so}_{4nm}$ 
be a nilpotent element with partition $((2m)^{2n})$.  
Then $\tilde\g^e\simeq\mathfrak L_{2m}(\spn,\sfr(\varpi_2)\oplus \odin)\simeq 
\mathfrak L_{2m}(\spn,\sfr(\varpi_2))\dotplus \z_m$.
Here $\spn\oplus(\sfr(\varpi_2)\oplus \odin)=\mathfrak{gl}_{2n}$ and $\tilde\g^e$ is a contraction of $m(\mathfrak{gl}_{2n})$.
\item[4$^o$.]  Let  $e\in \tilde\g=\mathfrak{sp}_{2n(2m+1)}$ 
be a nilpotent element with partition $((2m+1)^{2n})$.  
Then $\tilde\g^e\simeq\mathfrak L_{2m+1}(\spn,\sfr(\varpi_2)\oplus \odin)\simeq 
\mathfrak L_{2m+1}(\spn,\sfr(\varpi_2))\dotplus \z_m$.
Here $m(\mathfrak{gl}_{2n})\dotplus\spn\leadsto \tilde\g^e$.
\item[5$^o$.] Let $e\in\tilde\g= \mathfrak{so}_{2n+2}$ be a nilpotent element with partition $(3,1^{2n-1})$.
Then $\tilde\g^e\simeq (\mathfrak{so}_{2n-1}\ltimes\sfr(\varpi_1))\dotplus\bbk e$, where the
first summand is a contraction of $\sone$.
\item[6$^o$.] Let $e\in \tilde\g=\GR{E}{7}$ be a nilpotent element with weighted 
Dynkin diagram\\
$\left(\text{\begin{E7}{2}{0}{0}{0}{0}{0}{0}
\end{E7}}\right)$. Then $\tilde\g^e=(\GR{F}{4}\ltimes\sfr(\varpi_1))\dotplus
\bbk e$.
\item[7$^o$.]  For $n=6,7,8$, let $e\in \tilde\g=\GR{E}{n}$ be nilpotent elements with 
weighted Dynkin diagrams
$\left(\text{\begin{E6}{2}{0}{0}{0}{2}{0}\end{E6}}\right), 
\left(\text{\begin{E7}{2}{2}{0}{0}{0}{2}{0}\end{E7}}\right), 
\left(\text{\begin{E8}{2}{2}{2}{0}{0}{0}{2}{0}\end{E8}}\right)$,
respectively. Then $\tilde\g^e\simeq 
(\GR{G}{2}{\ltimes}\sfr(\varpi_1){\ltimes}\sfr(\varpi_1))\dotplus\z_{n-4}$, where the 
first summand is a $\BZ_3$-contraction of $\GR{D}{4}$.
\end{itemize}
Theorem~\ref{teor1+2} applies  to items $3^o$--$7^o$, and to $1^o$--$2^o$ if $n$ is odd.
The centralisers in $2^o$ and $4^o$ are quadratic.
\end{ex}

\begin{rmk}
If $|\theta|=2$  and $n=1$,  then the algebra 
$\tinv(\g_0\ltimes\g_1,\ad)$ is always polynomial and the basic invariants can explicitly be
described~\cite[Theorem\,6.2]{p05}. However, the equality 
$\tinv(\g_0\ltimes\g_1,\ad)=\gr_\bullet (\bbk[\g]^\g)$ holds if and only if 
$\g_0\cap\co^{reg}\ne\varnothing$.
\end{rmk}

\section{Good cases for coadjoint representations} 
\label{sect:coadj}

\noindent
In this section, we prove a counterpart of Theorem~\ref{teor1+2} for invariants of 
coadjoint representations of $\gp{nk}_0$. To this end, we need some preparations.

The {\it index\/} of a Lie algebra $\q$, $\ind\q$, is the minimal dimension of stabilisers 
of elements of $\q^*$ with respect to the coadjoint representation.
Let $\q^*_{reg}$ be the set of regular elements of $\q^*$, i.e., those with minimal
dimension of the stabiliser. We say that $\q$ has the {\it codim--2} (resp.
{\it codim--3}) property,
if $\codim(\q^*\setminus \q^*_{reg})\ge 2$  (resp. $\ge 3$).
Set $b(\q)=(\dim\q+\ind\q)/2$.
We will need the following result, which is explicitly stated in \cite[Thm.\,1.2]{coadj} and
based on an earlier work of Odesskii-Rubtsov \cite{or}.

\begin{thm}    \label{thm:coadj}
Suppose $\q$ has the codim--2 property and $\trdeg \cs(\q)^\q=\ind\q$. 
Set  $l=\ind\q$. Let $f_1,\dots,f_l\in \cs(\q)^\q$ be arbitrary
homogeneous algebraically independent polynomials. Then 
\begin{itemize}
\item[\sf (i)] \ $\sum_{i=1}^l\deg f_i \ge b(\q)$;
\item[\sf (ii)] \  If\/  $\sum_{i=1}^l\deg f_i = b(\q)$, then
$\cs(\q)^\q$ is freely generated by $f_1,\dots,f_l$ and
$\xi\in\q^*_{reg}$ if and only if $(\textsl{d}f_1)_\xi,\dots,(\textsl{d}f_l)_\xi$ are linearly
independent.
\end{itemize}
\end{thm}

\noindent In order to apply this result to algebras $\gp{nk}_0$, we must have the
{\it codim--2} property and know the index of $\gp{nk}_0$.

Given a periodic $\theta\in{\rm Aut}(\g)$, we have the flat quotient morphism 
$\pi:\g_1\to \g_1\md G_0$ (see Introduction). 
The fibre  $\fN_1:=\pi^{-1}(\pi(0))$ consists of all nilpotent 
elements in $\g_1$. 
It can be reducible; moreover, if $|\theta|\ge 3$, then some components can be 
reduced, while some other  not.

Let us say that $\theta$ (or the corresponding $\BZ_k$-grading) is
\begin{itemize}
\item  $\eus S$-{\it regular}, if $\g_1$ contains  regular semisimple elements of $\g$;
\item  $\eus N$-{\it regular}, if $\g_1$ contains  regular nilpotent elements of $\g$ (i.e.,
$\g_1\cap\co^{reg}\ne\varnothing$);
\item  {\it very\/} $\eus N$-{\it regular}, if $\co^{reg}$ meets each irreducible component of\/ $\fN_1$.
\end{itemize}
Some structure results on $\eus S$-  and $\eus N$-{regular}  gradings are obtained  in \cite{theta}.

\begin{lm}    \label{lm:big-open}
Suppose   $\theta$ is both $\eus S$-regular and very $\eus N$-regular.
Then  $\codim(\g_1\setminus (\g_1\cap \g_{reg}))\ge 2$. 
\end{lm}\begin{proof}
Let $\pi^{-1}(\beta)=\eus F_\beta$ be an arbitrary fibre of $\pi$  
and $\co_\beta$  the dense $G_0$-orbit in an irreducible component of
$\eus F_\beta$. Since $\dim\eus F_\beta=\dim\fN_1$,
the {\it associated cone\/} of $\co_\beta$ (see \cite[\S\,3]{bk})
is an irreducible component of $\fN_1$. Then the assumption on $\fN_1$ 
implies that $\co_\beta\subset \g_{reg}$. In other words, if $\co\subset\g_1$ is a $G_0$-orbit
of maximal dimension, then $\co\subset\g_{reg}$. 
The assumption on regular semisimple elements shows that a generic fibre of $\pi$ is a 
(closed) $G_0$-orbit, i.e., the whole such fibre belongs to 
$\g_1\cap\g_{reg}$. The union of all other fibres is
a proper subvariety of $\g_1$ (actually, it is a divisor). However, the open
$G_0$-orbits in all other fibres belong to $\g_{reg}$, too. Hence the complement is
of codimension at least 2, as required.
\end{proof}

\begin{prop}    \label{prop:codim-2}
{\sf (i)} \ Suppose that\/  $\g_1\cap\g_{reg}\ne \varnothing$. 
Then $\ind\gp{nk}_0=n{\cdot}\rk\g$.

\noindent
{\sf (ii)} \ If\/ $\theta$ is both $\eus S$-regular and very $\eus N$-regular, then 
$\gp{nk}_0$ has the {\it codim--2} property.
\end{prop}
\begin{proof}  (i) 
Set $\q=\gp{nk}_0$. Then $\q^*=\gp{nk}_1$. Recall that
\[
\begin{split}
\gp{nk}_0=& \,\g_0\ltimes\g_1\ltimes\ldots \ltimes \g_{k-1} \\
\gp{nk}_1=& \,\g_1\ltimes\g_2\ltimes\ldots \ltimes \g_0  \quad  \text{ ($nk$ factors in 
both cases)}
\end{split}
\]
Take any  $x\in \g_1\cap\g_{reg}$. Write $\g_i^x$ (resp. $\g^x$) for the centraliser of
$x$ in $\g_i$ (resp. $\g$). Then $\g^x=\oplus_{i=0}^{k-1}\g_i^x$ and $\dim\g^x=\rk\g$.
Consider $\boldsymbol\xi=(x,0,\dots,0)$ as element of $\q^*$. Then the stabiliser of 
$\boldsymbol\xi$ in $\q$ is $\g_0^x\ltimes\g_1^x\ltimes\dots \ltimes \g_{k-1}^x$, i.e.,
we get $\g_i^x$ inside every component $\g_i$ occurring in $\q$ . Since each 
$\g_i$ occurs $n$ times, $\dim \q^{\boldsymbol\xi}=
n{\cdot}\rk\g$. Therefore $\ind\q\le n{\cdot}\rk\g$. On the other hand,
$\q$ is a  contraction of $n\g$ (Theorem~\ref{thm:main-contr}(i)). Hence
$\ind\q\ge \ind (n\g)=n{\cdot}\rk\g$. 
Notice that we also proved that if $x\in\g_1\cap\g_{reg}$, then
$\boldsymbol\xi=(x,0,\dots,0)\in \q^*_{reg}$.

(ii) Set $\g_1^{reg}=\g_1\cap\g_{reg}$. By Lemma~\ref{lm:big-open}, we have 
$\codim (\g_1\setminus \g_1^{reg})\ge 2$.
Now, let $\boldsymbol\xi=(\xi_1,\xi_2,\dots,\xi_{nk})\in \q^*$, where $\xi_i\in\g_{\bar i}$. 
We claim that  if
$\xi_1\in \g_1^{reg}$, then $\boldsymbol\xi\in\q^*_{reg}$. This yields the desired 
{\it codim--2} property. Hence it suffices to prove the claim.
Consider $\boldsymbol\xi(t)=(\xi_1,t\xi_2,\dots,t^{nk-1}\xi_{nk})$, 
where  $t\in\bbk$.
It is easily seen that if $(x_0,x_1,\dots,x_{nk-1})\in \q^{\boldsymbol\xi(t)}$ for $t\ne 0$, then
$(t^{nk-1}x_0,t^{nk-2}x_1,\dots,x_{nk-1})\in \q^{\boldsymbol\xi}$. 
It follows that, for $t\ne 0$, $\dim\q^{\boldsymbol\xi(t)}$
does not depend on $t$. Because $\lim_{t\to 0}{\boldsymbol\xi(t)}=
(\xi_1,0,\dots,0)\in \q^*_{reg}$, 
we conclude that
all elements $\boldsymbol\xi(t)$ are regular. That is, 
$(\xi_1,\xi_2,\dots,\xi_{nk})\in\q^*_{reg}$ whenever  $\xi_1\in \g_1^{reg}$. 
\end{proof}

\begin{ex}   \label{ex:bad}
It can happen that  $\theta$ is  both $\eus S$- and  $\eus N$-regular,
but not very $\eus N$-regular.
This may lead to the absence of the {\it codim--2} property for $\gp{|\theta|}_0$.
Consider $\g=\mathfrak{sp}_4$ and an automorphism of
order 4 such that $\g_0$ is a Cartan subalgebra, say $\te$.
If $\ap,\beta$ are simple roots with respect to $\te$ ($\ap$ is short), so that $2\ap+\beta$
is the highest root, 
then $\g_1$ is the sum of root spaces corresponding to $\ap,\beta$, and $-2\ap{-}\beta$.
Therefore  
$\g_1$ contains regular semisimple and regular nilpotent elements of $\g$. However,
a direct verification shows that the corresponding $\BZ_4$-contraction does not
have the {\it codim--2} property. In this case, $\fN_1$ has three
irreducible components and $\co^{reg}$ meets only two of them.
\end{ex}

\begin{thm}       \label{teor3}
Suppose $\theta$ is both $\eus S$-regular and very $\eus N$-regular and $n\in\BN$.
Then 
 $\tinv(\gp{nk}_0, \ads)$ is a polynomial algebra of Krull dimension $n{\cdot}\rk\g$.
Moreover, $\eus L^{\bullet}(\cs(n\g)^{n\g})=\tinv(\gp{nk}_0, \ads)$.
\end{thm}
\begin{proof}  
First, assume that $n=1$. 
 In view of Theorem~\ref{thm:coadj} and Proposition~\ref{prop:codim-2},
it suffices to show that there are basic invariants $F^{(1)},\dots,F^{(\rk\g)}\in \cs(\g)^{\g}$ 
such that $\{\gr^\bullet F^{(i)}\}$ remain algebraically independent. 

Take the basic invariants in $\cs(\g)^{\g}$ such that each $F^{(i)}$ is an eigenvector of 
$\theta$. Let $\deg F^{(i)}=d_i$.
We want to better understand the structure of $\gr^\bullet F^{(i)}$.
Take $e\in\co^{reg}\cap\g_1$. 
By \cite{ko63}, the differentials
$\{\textsl{d}F^{(i)}_e\}$ are linearly independent. In particular, $\textsl{d}F^{(i)}_e\ne 0$
for each $i$. Consider the  polygrading of $\cs(\g)$ corresponding
to the decomposition $\g=\bigoplus_{i=0}^{k-1} \g_i$. 
We also regard $\g_1$ as the first 
factor of $\gp{k}_1\simeq \gp{k}_0^*$. Since $\g_{k-1}^*\simeq \g_1$
and $\g_{k-1}$ is the last factor of $\gp{k}_0$, the condition
$\textsl{d}F^{(i)}_e\ne 0$ implies that $F^{(i)}$ has a nonzero summand of the form 
\\[1ex]
\hbox to \textwidth{ \ ({\color{blue}$\clubsuit$}) \hfil 
$F^{(i)}_{(0,\dots,0,d_i)}$\quad or \quad  $F^{(i)}_{(\dots,1,\dots,d_i-1)}$. \hfil}

\noindent
We have $k-1$ possibilities for  the position of `1' in the second expression, hence
totally $k$ possibilities  in ({\color{blue}$\clubsuit$}). In fact, there is the following
precise assertion:

\begin{lm}   \label{lm:precise-gr}
(i) A homogeneous $\theta$-eigenvector $F\in\cs(\g)$ can have at most one nonzero 
summand of the form ({\color{blue}$\clubsuit$});
(ii) If this is the case, then $\gr^\bullet F$ contains that summand.
\end{lm}\begin{proof}
(i) This follows from the fact that, for the nonzero summnads $F_{(i_0,\dots,i_{k-1})}$, 
$i_0+\ldots +i_{k-1}=\deg F$  and all the sums $\sum_{j=1}^{k-1}j i_j$ have one and 
the same residue $\pmod k$, which is determined by the $\theta$-eigenvalue of $F$.
The $k$ possibilities in ({\color{blue}$\clubsuit$}) just correspond to $k$ possible eigenvalues.

(ii)  If $\deg F=d$ and $F$ has the summand $F_{(0,\dots,0,d)}$, then the latter is clearly
$\gr^\bullet F$. If $F$ has the summand $F_{(\dots,1,\dots,d-1)}$ (with `1' in position $s$, 
$0\le s\le k-2$), then 
the summands $F_{(i_0,\dots,i_{k-1})}$ occurring in $\gr^\bullet F$ 
should satisfy the relations
\[
    \begin{cases}   \sum_{j=0}^{k-1} i_j =d , \\
                             \sum_{j=1}^{k-1} j i_j \equiv s +(d-1)(k-1) \pmod k , \\
                    \sum_{j=1}^{k-1} j i_j \ \text{  is maximal posiible}.
   \end{cases}
\]
It is not hard to prove that the maximal value of the last sum is $s +(d-1)(k-1)$. Hence
one of the solutions is  $(\dots,\underset{s}{1},\dots,d-1)$.
\end{proof}

\noindent It follows from Lemma~\ref{lm:precise-gr}(ii)  
that $\textsl{d}F^{(i)}_e=\textsl{d}(\gr^\bullet F^{(i)})_e$.
Hence $\{\gr^\bullet F^{(i)}\}$ remain algebraically independent. This completes the proof 
of theorem for $n=1$.

For arbitrary $n\in\BN$, we consider $\widetilde\theta\in {\rm Aut}(n\g)$ and the 
eigenspaces $(n\g)_j$, which are described in Lemma~\ref{lm:ext-theta}. Here
$(n\g)_1=\{(x,\mu x,\dots, \mu^{n-1}x)\mid x\in \g_1\}$, where $\mu=\sqrt[nk]1$, and we 
work with 
$\boldsymbol e=(e,\mu e,\dots, \mu^{n-1}e)$, which is a regular nilpotent element of $n\g$.
We regard $(n\g)_1$ as the first factor in $\gp{nk}_1$ and choose  basic invariants 
in $\cs(n\g)^{n\g}$ that are $\widetilde\theta$-eigenvectors. The rest is the same.
\end{proof}

\noindent
The property of being ``very $\eus N$-regular'' is difficult to check directly. 
There are, however, useful sufficient conditions.

For $\vert \theta\vert=2$, all irreducible component of $\fN_1$ are conjugate with respect 
to the action of certain (non-connected) group containing $G_0$, see~\cite[Theorem\,6]{kr71}. 
Therefore $\eus N$-regularity coincides with very $\eus N$-regularity.
Furthermore, an involution is $\eus S$-regular if and only if it is $\eus N$-regular~\cite{an}. 
Thus, it suffices to assume that $\g_1\cap\co^{reg}\ne \varnothing$. Finally,
an involution has the last property if and only if the corresponding Satake diagram 
has no black nodes.

To state another sufficient condition, we recall that, for $\g$  simple, 
the set of basic invariants in $\bbk[\g]^G$ contains a unique polynomial of maximal degree.
This degree equals the Coxeter number of $\g$, denoted $\mathsf c(\g)$.
Let $F_{\mathsf c(\g)}$ be such a basic invariant.
It is known that $\textsl{d}(F_{\mathsf c(\g)})_e=0$ for any $e\in \fN\setminus \co^{reg}$
(see a description of the ideal of $\fN\setminus \co^{reg}$ in \cite[4.7--4.9]{broer}).

\begin{prop}   \label{pr:suff}
Let $\g$ be simple. Suppose $\theta$ is $\eus N$-regular 
and $F_{\mathsf c(\g)}\vert_{\g_1}\not\equiv 0$. If either {\rm (i)} \ $\g_0$ is semisimple or
{\rm (ii)} \ $G_0\subset SL(\g_1)$ and $(G_0:\g_1)$ is locally free, 
then $\theta$ is very $\eus N$-regular.
\end{prop}\begin{proof}
Since $\g_1\cap\co^{reg}\ne \varnothing$, the restriction homomorphism
$\bbk[\g]^G \to \bbk[\g_1]^{G_0}$ is onto, and $F_{\mathsf c(\g)}\vert_{\g_1}$ is a basic invariant
in  $\bbk[\g_1]^{G_0}$ \cite[Theorem\,3.5]{theta}.
Let $\co$ be a dense $G_0$-orbit in an irreducible component of $\fN_1$.
In other words, $\co$ is a nilpotent $G_0$-orbit of maximal dimension.
Then in both cases, the differentials of basic
invariants in $\bbk[\g_1]^{G_0}$ are linearly independent at any $v\in \co$. For (i)
(resp. (ii)), we refer to
\cite[Cor.\,5(i)]{p85} (resp. \cite[Cor.\,1]{p84}).
In particular, we have $\textit{d}(F_{\mathsf c(\g)})_v\ne 0$. Hence $v\in \co^{reg}$.
\end{proof}
\begin{rem}
In Example~\ref{ex:bad}, the condition that $G_0\subset SL(\g_1)$ is not satisfied.
\end{rem}
Here are some other related results:
\begin{prop}    \label{pr:other_use}
{\ }\phantom{aaa} \\
{\sf (i)} \ If $\theta$ is $\eus S$-regular  and the $G_0$-action on $\g_1$ is locally free, then
$\theta$ is $\eus N$-regular \cite[Thm.\,4.2(iii)]{theta}; 
{\sf (ii)} \
If $\theta\in {\rm Int}(\g)$ and $\theta$ is $\eus N$-regular, then 
$F_{\mathsf c(\g)}\vert_{\g_1}\not\equiv 0$ if and only if\/  
$|\theta|$ divides $\mathsf c(\g)$ \ \cite[Cor.\,3.6]{theta}.
\end{prop}

In the rest of the section, we provide  examples of $\eus S$-regular and
very $\eus N$-regular automorphisms and thereby examples where Theorem~\ref{teor3}
applies.

\begin{ex} We give some serial examples of $\eus S$-regular and very $\eus N$-regular periodic automorphisms related to classical algebras.

1)  The Lie algebra $\mathfrak{gl}_{nk}$ has an automorphism of order $k$
such that $\g_0=k\gln$. In more invariant terms, let 
$V_i$ be a $n$-dimensional space, $i=1,\dots,k$, and 
$\g=\mathfrak{gl}(V_1\oplus\dots\oplus V_k)$. Define $A\in
GL(V_1\oplus\dots\oplus V_k)$ by $A\vert_{V_i}=\zeta^i{\cdot}{\rm id}$.
Let $\theta$ be the inner automorphism of $\g$ determined by $A$.
Then $\g_0=\mathfrak{gl}(V_1)\dotplus\ldots\dotplus \mathfrak{gl}(V_k)$ and 
$\g_1= \bigoplus_{i=1}^k {\rm Hom}(V_i,V_{i+1})$, where $V_{k+1}=V_1$.
Here $\dim\g_1=kn^2$. The generic stabiliser in $\g_0$ is $\te_n$ (the diagonally
embedded Cartan subalgebra). Using the matrix realisation, one easily verifies
that $\theta$ is both $\eus S$- and $\eus N$-regular. A more careful argument shows
that $\fN_1$ has $k$ irreducible components, and each contains regular nilpotent
elements of $\g$. Therefore $\theta$ is very $\eus N$-regular.

2)  The algebra $\GR{D}{4m+3}$ has an inner automorphism $\theta$ of order 4 such that
$G_0=\GR{D}{m+1}\times\GR{D}{m+1}\times\GR{A}{2m}\times T_1$.
The corresponding Kac's diagram is 

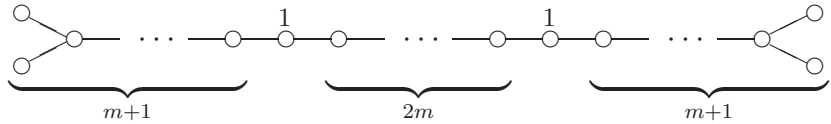
\begin{center}
\begin{picture}(350,45)(0,-15)
\put(30,8){\circle{6}}
\put(290,8){\circle{6}}
\multiput(90,8)(20,0){3}{\circle{6}}
\multiput(190,8)(20,0){3}{\circle{6}}
\multiput(10,-2)(0,20){2}{\circle{6}}
\multiput(310,-2)(0,20){2}{\circle{6}}
\multiput(33,8)(240,0){2}{\line(1,0){14}}
\multiput(73,8)(20,0){4}{\line(1,0){14}}
\multiput(173,8)(20,0){4}{\line(1,0){14}}
\put(13,0){\line(2,1){14}}
\put(13,16){\line(2,-1){14}}

\put(293,10){\line(2,1){14}}
\put(293,6){\line(2,-1){14}}

\put(54,5){$\cdots$}
\put(154,5){$\cdots$}
\put(254,5){$\cdots$}
\put(107,13){\footnotesize $1$}
\put(207,13){\footnotesize $1$}
\put(5,-7){$\underbrace%
{\mbox{\hspace{90\unitlength}}}_{m+1}$}
\put(125,-7){$\underbrace%
{\mbox{\hspace{70\unitlength}}}_{2m}$}
\put(225,-7){$\underbrace%
{\mbox{\hspace{90\unitlength}}}_{m+1}$}
\end{picture}
\end{center}
\noindent
We refer to  \cite[\S\,8, Prop.\,17]{vi76} for a complete account on Kac's diagrams of
periodic automorphisms. (Partial explanations are also given in \cite[Example\,4.5]{theta},
where we have drawn black nodes in place of  nodes with labels `$1$'.)
Here $\dim\g_0=(2m+1)(4m+3)$ and $\dim\g_1=(2m+1)(4m+4)$. Let us prove that
$\theta$ is $\eus S$-regular and very $\eus N$-regular.

Since $\mathsf c(\GR{D}{4m+3})=8m+4$, $|\theta|$ divides it.
The representation of $G_0$ in $\g_1$ can be read off the Kac's diagram. 
Here $\g_1$ is the sum of two simple $G_0$-modules of the {\sl same\/} dimension, and
$T_1$ acts with opposite weights on them. Therefore $G_0\subset SL(\g_1)$.
One also readily verifies that 
the action $(G_0:\g_1)$ is stable and locally free.
Therefore $\dim\ce=\dim\g_1-\dim\g_0=2m+1$.  
For stable $\theta$-groups, the dimension
of a generic semisimple $G$-orbit meeting $\g_1$ can be computed as
$|\theta|(\dim\g_1-\dim\ce)$ (cf. \cite[Prop.\,2.1(i)]{theta}), 
which is equal in this case to $\dim\g -\rk\g$. Hence, $\theta$ is $\eus S$-regular.
Combining Proposition~\ref{pr:suff}(ii) and \ref{pr:other_use}, we then conclude that
$\theta$ is very $\eus N$-regular.
\end{ex}

\begin{ex} 
The following table contains some sporadic examples, mostly
for exceptional Lie algebras. Here $\g_0$ is always semisimple and
$|\theta|$ divides $\mathsf c(\g)$; 
$\ind(\theta)$ denotes the order of $\theta$ in ${\rm Aut}(\g)/{\rm Int}(\g)$. In the last two 
columns, the dimension of Cartan
subspaces and the generic stabiliser for the $\g_0$-module $\g_1$ are displayed
(this information is borrowed from \cite[\S\,9]{vi76}). The $G_0$-action on $\g_1$ is locally
free if and only if this generic stabiliser is trivial.

\noindent
In most cases, the proof is similar to that given in the previous example. The only 
exception is the second item for $\GR{E}{7}$, where $\theta$ is not an involution 
and $(G_0:\g_1)$ is not locally free. However, in this case $\fN_1$ appears to be 
irreducible \cite[\S\,4]{gati}, hence $\eus N$-regularity is sufficient.

\begin{longtable}{ccccccc}
$\g$ & $|\theta|$  & $\ind(\theta)$ & $\g_0$  & $\g_1$ & $\dim\ce$ & $\es$ \\  \hline\hline 
\nopagebreak \samepage
$\GR{E}{6}$  & 2  & 1 &  $\GR{A}{5}\times\GR{A}{1}$  & $\sfr(\varpi_3)\otimes\sfr(\varpi_1)$ 
& 4 & $\te_2$ \\
  & 3  & 1 &  $\GR{A}{2}{\times} \GR{A}{2}{\times} \GR{A}{2}$  
& $\sfr(\varpi_1){\otimes}\sfr(\varpi_1){\otimes}\sfr(\varpi_1)$ & 3 & $\{0\}$ \\
  & 2  & 2 &  $\GR{C}{4}$  & $\sfr(\varpi_4)$ & 6 & $\{0\}$ \\
  & 4  & 2 &  $\GR{A}{3}\times\GR{A}{1}$  & $\sfr(2\varpi_1)
\otimes\sfr(\varpi_1)$ & 2 & $\{0\}$ \\  \hline
$\GR{E}{7}$  & 2  & 1 &  $\GR{A}{7}$  & $\sfr(\varpi_4)$ & 7 & $\{0\}$ \\
  & 3  & 1 &  $\GR{A}{5}\times \GR{A}{2}$  & $\sfr(\varpi_2)\otimes\sfr(\varpi_1)$ & 3 
& $\te_1$  \\ \hline
$\GR{E}{8}$  & 2  & 1 &  $\GR{D}{8}$  & $\sfr(\varpi_7)$ & 8 & $\{0\}$ \\
      & 3  & 1 &  $\GR{A}{8}$  & $\sfr(\varpi_3)$ & 4 & $\{0\}$ \\
  & 5  & 1 &  $\GR{A}{4}\times \GR{A}{4}$  & $\sfr(\varpi_2)\otimes\sfr(\varpi_1)$ 
& 2 & $\{0\}$ \\ \hline
$\GR{F}{4}$  & 2  & 1 &  $\GR{C}{3}\times\GR{A}{1}$  & $\sfr(\varpi_3)
\otimes\sfr(\varpi_1)$  & 4 & $\{0\}$ \\
       & 3  & 1 &  $\GR{A}{2}\times\GR{A}{2}$  & $\sfr(2\varpi_1) \otimes\sfr(\varpi_1)$ 
& 2 & $\{0\}$ \\ \hline
$\GR{G}{2}$  & 2  & 1 &  $\GR{A}{1}\times\GR{A}{1}$  & $\sfr(3\varpi_1)
\otimes\sfr(\varpi_1)$  & 2 & $\{0\}$ \\ 
$\GR{D}{4}$  & 3  & 3 &  $\GR{A}{2}$  & $\sfr(3\varpi_1)$ & 2 & $\{0\}$ \\  \hline
\end{longtable} 
\end{ex} 

\section{Some remarks and open problems} 
\label{sect:end}

\subsection{Invariants of the nilpotent radical}
In \cite{p05}, we explore several instances of representations of non-reductive Lie algebras 
$\q$ such that, for certain nilpotent ideal $\rr\triangleleft\q$, the invariants of $\rr$ are
polynomial. This includes the adjoint and coadjoint representations of $\q$, i.e., the 
algebras  $\bbk[\q]^\rr$ or $\cs(\q)^\rr$.  
Often, $\rr=\mathfrak R_u(\q)$. This fact was used as a step
toward proving the polynomiality of  algebras of $\q$-invariants.
(See Theorems 6.2, 7.1, and 11.1 in loc. cit.) 

Similar approach applies to our results in Sections~\ref{sect:adj} and \ref{sect:coadj}.
For instance, in order to describe $\tinv(\gp{k}_0,\ad)$, one can first consider the algebra of
$\mathfrak R_u(\gp{k}_0)$-invariants in $\bbk[\gp{k}_0]$, say $\ca$.
Under the assumption of Theorem~\ref{teor1+2}, it can be shown that 
$\ca$ is polynomial, of Krull dimension 
$\dim\g_0 +(\rk\g-\rk\g_0)$. Furthermore, the induced representation of $\g_0\simeq
\gp{k}_0/\mathfrak R_u(\gp{k}_0)$ in
$\spe(\ca)$ is isomorphic to the adjoint representation of $\g_0$ plus the trivial
representation of dimension $\rk\g-\rk\g_0$. In this way, one obtains another proof
of polynomiality of $\tinv(\gp{k}_0,\ad)$.
The reason for success in this and other similar cases is that one can explicitly
construct a natural  set of elements of $\ca$ (presumable basic invariants). Then
using Igusa's lemma (see e.g. \cite[Lemma\,6.1]{p05}), one proves that these elements are algebraically independent and generate the algebra $\ca$. The details will appear elsewhere.

\subsection{Index of fixed point subalgebras}
Let us summarise what is known about the index of algebras $\gm_0$. Recall that 
$|\theta|=k$.

\noindent
Since $\gp{nk+1}_0$ is quadratic, the description of 
$\tinv(\gp{nk+1}_0,\ad)$  in Theorem~\ref{teor1+2}(ii) shows that
$\ins(\gp{nk+1}_0)=n{\cdot}\rk\g+\rk\g_0$. That is, index does not change under the
contraction $n\g\dotplus\g_0\leadsto \gp{nk+1}_0$ whenever 
$\g_0\cap\co^{reg}\ne\varnothing$. (Actually, this equality can be proved if 
$\g_0\cap\g_{reg}\ne\varnothing$.)
Similarly, if $\g_1\cap\g_{reg}\ne\varnothing$, then $\ins(\gp{nk}_0)=n{\cdot}\rk\g$
(Proposition~\ref{prop:codim-2}).

If $|\theta|=2$, then we always have $\ins(\g_0\ltimes\g_1)=\ind\g$
\cite[Corollary\,9.4]{p05}. The reason is that here $G/G_0$ is a spherical homogeneous 
space.
From this, it is not hard to deduce that
$\ins(\mathfrak L_{2n}(\g_0,\g_1))=n{\cdot}\rk\g$ and
$\ins(\mathfrak L_{2n+1}(\g_0,\g_1))=n{\cdot}\rk\g+\rk\g_0$.
(See notation introduced in Example~\ref{real-life}.)
However, for $|\theta|\ge 3$,  the general answer is not known. 

\begin{prob}
Compute $\ind(\gp{nk}_0)$ and  $\ind(\gp{nk+1}_0)$ for an
arbitrary $\theta\in{\rm Aut}(\g)$ with $|\theta|=k\ge 3$. Is it true that\/ 
$\ins(\gp{nk+1}_0)=n{\cdot}\rk\g+\rk\g_0$ and\/ $\ins(\gp{nk}_0)=n{\cdot}\rk\g$?
\end{prob}
\noindent
(The existence of contractions $n\g\dotplus\g_0\leadsto \gp{nk+1}_0$ and 
$n\g\leadsto \gp{nk}_0$ shows that in both cases the inequality ``$\ge$" holds.)
For instance, consider the outer automorphism of $\GR{D}{4}$ of order 3
whose fixed point subalgebra is $\GR{G}{2}$. The corresponding 
$\BZ_3$-contraction is $\q=\GR{G}{2}{\ltimes}\sfr(\varpi_1){\ltimes}\sfr(\varpi_1)$.
It occurs in Example~\ref{real-life}(7$^o)$. 
Here $\g_1\simeq \sfr(\varpi_1)$ does not contain regular elements of $\GR{D}{4}$.
However, it is not hard to verify that $\ind\q=4$.

\subsection{On Poisson-commutative subalgebras}
A subalgebra $\ca$ of $\cs(\q)$ is said to be  {\it Poisson-commutative\/} if $\{f,g\}=0$
for all $f,g\in \ca$. There is a procedure 
(the so-called {\it argument shift method}, see e.g.~\cite{codim3}) 
for constructing ``large'' Poisson-commutative 
subalgebras of $\cs(\q)$, which begins with $\cs(\q)^\q$ and a $\xi\in\q^*_{reg}$.
The resulting subalgebra is denoted by $\cf_\xi(\cs(\q)^\q)$.
It is proved in \cite{codim3} that if 
{\sf (1)} $\cs(\q)^\q$ is polynomial, 
{\sf (2)} the sum of degrees of the basic invariants of $\cs(\q)^\q$ equals $b(\q)$, and 
{\sf (3)} $\q$ has the {\it codim--3} property, then
$\cf_\xi(\cs(\q)^\q)$ is a maximal Poisson-commutative subalgebra for any 
$\xi\in\q^*_{reg}$. Furthermore, $\cf_\xi(\cs(\q)^\q)$ is a polynomial algebra of 
Krull dimension $b(\q)$.

Our goal is to realise when that result applies to algebras $\gm_0$.
First of all, Lie algebras $\q$ occurring in Theorems~\ref{teor1+2}(ii) and \ref{teor3}
satisfy properties {\sf (1)} and {\sf (2)}. We also proved that $\gp{nk}_0$ has 
the {\it codim--2} property (Proposition~\ref{prop:codim-2}(ii)). However, the
{\it codim--3} property does not always holds for $\gp{nk}_0$. But for algebras
$\gp{nk+1}_0$ the situation is better.

\begin{prop}
Suppose  $\theta$ has the property that $\g_0\cap \co^{reg}\ne\varnothing$.
Then $\gp{nk+1}_0$ has the {\it codim--3} property.
\end{prop}\begin{proof}
The proofs of Lemma~\ref{lm:big-open} and  
Prop.~\ref{prop:codim-2} can be adapted to this situation. Recall that 
\[
   \gp{nk+1}_0\simeq \gp{nk+1}_0^*=\g_0\ltimes\g_1\ltimes\ldots \ltimes \g_{k-1}
\ltimes \g_0  \quad  \text{ ($nk+1$ factors).}
\]
As in Lemma~\ref{lm:big-open}, we prove that if $x\in\g_0$ is regular in $\g_0$, then
it is regular in $\g$. (In doing so, we use the assumption $\co^{reg}\cap\g_0\ne
\varnothing$ and the fact that the nilpotent cone in $\g_0$ is irreducible.)
Then, as in Prop.~\ref{prop:codim-2}, we prove that 
$\boldsymbol\xi=(\xi_0,\xi_1,\dots,\xi_{nk})\in (\gp{nk+1}_0^*)_{reg}$
whenever $\xi_0\in(\g_0)_{reg}$. Since $\codim(\g_0\setminus (\g_0)_{reg})=3$,
we are done.
\end{proof}

\subsection{Flatness}    \label{subs:flat}
Although  we have found a number of periodic automorphisms of Takiff  algebras
such that $\tinv(\gm_0,\ad)$ or $\tinv(\gm_0,\ads)$ is polynomial, 
we do not have substantial results on the flatness of respective quotient morphisms. 
Actually, I believe that the  quotient morphisms are flat in the context of Theorems~\ref{teor1+2}
and \ref{teor3}.  Partial affirmative results for $\vert\theta\vert=2$ are contained in 
\cite[Theorem\,9.13]{p05} (the adjoint representation of $\g_0\ltimes\g_1$) 
and \cite[Sect.\,5]{coadj} (the coadjoint representation of $\g_0\ltimes\g_1$) .

\noindent
For the centraliser $\tilde\g^e$ from 
Example~\ref{real-life}(7$^o)$ and its adjoint representation,
we can also prove, using ad hoc methods, that the quotient morphism is flat.

\end{document}